\documentclass{amsart}

\usepackage{amssymb}
\usepackage{amsmath}
\usepackage{amsthm}
\usepackage{stmaryrd}
\usepackage{amsfonts,mathrsfs}
\title[Takayama extension]{A Takayama-type Extension Theorem}
\author{Dror Varolin}
\thanks{Partially supported by NSF grant DMS-0400909}
\address{
Department of Mathematics \newline \indent
Stony Brook University \newline \indent
Stony Brook, NY 11794}

\subjclass[2000]{32L10 14F10}

\newcommand{\noi}{\noindent}

\newcommand{\cc}{{\mathcal C}}

\newcommand{\co}{{\mathcal O}}

\newcommand{\si}{{\mathscr I}}

\newcommand{\vp}{\varphi} 
\newcommand{\ve}{\varepsilon}

\newcommand{\C}{{\mathbb C}}

\newcommand{\N}{{\mathbb N}}

\newcommand{\Q}{{\mathbb Q}}
\newcommand{\R}{{\mathbb R}}

\newcommand{\Z}{{\mathbb Z}}

\newcommand{\di}{\partial}
\newcommand{\dbar}{\bar \partial}

\newcommand{\relcomp}{\subset \subset}

\newcommand{\ii}{\sqrt{-1}}

\newcommand{\tensor}{\otimes}

\begin{document}
\maketitle

\theoremstyle{plain}
\newtheorem{thm}{\sc Theorem}
\newtheorem{lem}{\sc Lemma}[section]
\newtheorem{d-thm}[lem]{\sc Theorem}
\newtheorem{prop}[lem]{\sc Proposition}
\newtheorem{cor}[lem]{\sc Corollary}

\theoremstyle{definition}
\newtheorem{conj}[lem]{\sc Conjecture}
\newtheorem{defn}[lem]{\sc Definition}
\newtheorem{qn}[lem]{\sc Question}

\theoremstyle{definition}
\newtheorem{ex}{\sc Example}

\theoremstyle{remark}
\newtheorem*{rmk}{\sc Remark}
\newtheorem*{rmks}{\sc Remarks}
\newtheorem*{ack}{\sc Acknowledgment}

\section{Introduction}

Let $X$ be a compact complex algebraic manifold, $E \to X$ a holomorphic line bundle, and $Z \subset X$ a smooth codimension-1 submanifold.  The main goal of this paper is to establish sufficient conditions for extending sections of the pluri-adjoint bundles $m(K_Z+E|Z)$ from $Z$ to $X$.

Such an extension theorem was established by Takayama \cite[Theorem 4.1]{tak}, in the situation where $E$, seen as a $\Q$-divisor, is a sum of a big and nef $\Q$-divisor and a $\Q$-divisor whose restriction to $Z$ is Kawamata log terminal.  (The precise statement is Theorem \ref{tak-thm} below.)  In this paper, we wish to weaken the positivity hypotheses on $E$.

The following theorem is our main result.

\begin{thm}\label{main}
Let $X$ be a projective algebraic manifold, and $Z \subset X$ a smooth complex submanifold of codimension $1$.  Denote by $T \in H^0(X,Z)$ the canonical holomorphic section whose zero divisor is $Z$.  Let $E , B \to X$ holomorphic line bundles and assume there exist singular Hermitian metrics $e^{-\vp _Z}, \ e^{-\vp _E}$ and $e^{-\vp _B}$ for the line bundle associated to $Z$, for $E$ and for $B$ respectively, with the following properties:
\begin{enumerate}
\item[(R)]  The metrics $e^{-\vp _Z}$, $e^{-\vp _E}$ and $e^{-\vp _B}$ restrict to singular Hermitian metrics on $Z$.
\item[(B)]  The metric $e^{-\vp _Z}$ satisfies the uniform bound 
\[
\sup _{X} |T|^2e^{-\vp _Z} < +\infty.
\]
\item[(P)] 
There is an integer $\mu > 0$ such that 
\[
\ii \di \dbar ( \mu (\vp _E + \vp _B)) \ge \max (\ii \di \dbar \vp _Z , 0).
\]  
\item[(T)]  The multiplier ideal of $(\vp _Z+\vp _E)|Z$ is trivial:
$\si (e^{-(\vp _Z+\vp_E)}|Z) = \co _Z$.
\end{enumerate}
Then every global section of 
\[
H^0(Z, \co _Z (m(K_Z + E|Z)+B|Z )\tensor \si (e^{-(\vp _Z + \vp _E +\vp _B)}|Z))
\]
extends to a global holomorphic section in $H^0(X,m(K_X+Z+E)+B)$.
\end{thm}

\begin{rmk}
When we say that a section $s \in H^0(Z,K_Z+E|Z)$ extends to a section $S$ in $H^0(X,K_X+Z+E)$, we mean that 
\[
S|Z=s\wedge dT.
\]
\end{rmk}

\begin{rmk}
Recall that a singular Hermitian metric is a Hermitian metric for a holomorphic line bundle such that, if the metric is represented locally by $e^{-\vp}$, then ${\vp}$ is $L^1_{\ell oc}$.  In most situations, people deal with such metrics only when $\vp$ is plurisubharmonic.  Such $\vp$ are locally bounded above.  Here we have a metric in the picture that need not be plurisubharmonic, namely $e^{-\vp _Z}$.  Thus we add to our definition of singular Hermitian metrics the additional requirement that the local potentials $\vp$ be uniformly bounded above on their domain of definition.

Of course, with this convention, a singular metric for a holomorphic line bundle does not induce a singular metric for the dual bundle.  However, this asymmetry will not pose a problem for us.
\end{rmk}

An immediate corollary of Theorem \ref{main} is the following result.

\begin{cor}\label{no-mult}
Let the notation and hypotheses of Theorem \ref{main} hold.  In addition, assume that the metric $e^{-(\vp _Z + \vp _E+\vp _B)}|Z$ is locally integrable on $Z$.  Then the natural restriction map 
\[
H^0(X, m(K_X +Z+E)+B) \to H^0(Z,m(K_Z +E|Z)+B|Z)
\]
is surjective.  In particular, if $B = \co _X$ and $\vp _B \equiv 0$, then 
\[
H^0(X,m(K_X+Z+E)) \to H^0(Z,m(K_Z+E|Z))
\]
is surjective.
\end{cor}

\noi In Section \ref{tak-section} we shall deduce from Corollary \ref{no-mult} Takayama's Theorem as well as a more general algebro-geometric extension theorem (Theorems \ref{tak-thm} and \ref{tak-2} respectively).

\medskip

A typical way to extend sections on line bundles in algebraic geometry is through the use of vanishing theorems.  In the case of Takayama's Theorem, it is the Nadel Vanishing Theorem that provides a key step in the proof.  In order to use Nadel's Theorem, Takayama must assume that the divisor $E$ above is big and nef.  From the analytic perspective, Takayama's method requires that $E$ support a singular Hermitian metric having strictly positive curvature current.

The present paper came about because in complex analysis there is technology available for extension theorems that does not require strict positivity.  This technology arose first in the work of Ohsawa and Takegoshi, in the context of extension of holomorphic functions, square integrable with respect to a plurisubharmonic weight, from hyperplanes to a pseudoconvex domain in $\C ^n$ \cite{ot-87}.  Since that time there have been several extensions of the result.  The simplest general exposition was finally given by Siu \cite{s-02}, and used as one of two fundamental tools in establishing the celebrated deformation invariance of plurigenera.  Extensions to more singular settings were established by the author and McNeal in \cite{mv}.  

The positivity assumptions in Theorem \ref{main} are in some sense minimal, and in particular are insufficient to support the use of vanishing theorems.  A technique for dealing with such a situation was initiated by Siu in \cite{s-02}, and we are going to use a similar approach.  

As already suggested, at the foundations of such an approach lies an $L^2$ extension theorem of twisted canonical sections from the divisor $Z$ to the ambient manifold $X$.   In regard to extension of pluricanonical sections, so far the $L^2$-extension theorems that have been used are tailored to the case where the normal bundle of $Z$ is trivial.  In this setting, special properties of the canonical bundle allow one to obtain an extension theorem that require only the non-negativity of the twisting line bundle.  

In this paper we treat the situation in which the normal bundle of $Z$ is not trivial.  In view of the adjunction formula, which is a formula for local extension of canonical sections, we are forced to deal with the possible negative contribution from the curvature of the conormal bundle of the divisor $Z$.  While there are results in the literature (see, e.g., \cite{m,d}) that handle the case in which the normal bundle of $Z$ is non-trivial, these results are not nicely compatible with the condition (P) of Theorem \ref{main}.  However, new methods developed by the author and J. McNeal can be adapted to the situation of non-trivial normal bundle so as to handle a condition like (P).  

Let us be slightly more precise.  Suppose we are attempting to extend an $H$-twisted canonical section $s$ on $Z$ such that 
\[
\int _Z |s|^2 e^{-\kappa} < +\infty \quad \text{and} \quad s \wedge dT \in \si (e^{-(\vp _Z + \kappa)}|Z)
\]
for some singular Hermitian metric $e^{-\kappa}$.  Here $T$ is the canonical holomorphic section for the line bundle associated to $Z$, whose zero divisor is $Z$.  The need for a metric $e^{-\vp _Z}$ and the multiplier ideal requirement arises from the adjunction formula.  It turns out that the technique we use requires the following curvature condition:
\begin{enumerate}
\item[(C)] For some positive integer $\mu$, $\mu \ii \di \dbar \kappa \ge \max (\ii \di \dbar \vp _Z , 0)$.
\end{enumerate}
Under such a condition, we establish in Section \ref{l2-section} an $L^2$ extension theorem, namely Theorem \ref{ot-thm}, which we consider to be the main new contribution of the present paper. 

Ultimately, the extension theorem is used in the following way.  Suppose we are given a singular metric $e^{-\vp _Z}$ and a global section $T$ of the line bundle associated to the divisor $Z$, metrics $e^{-\vp _E}$ and $e^{-\vp _B}$ for the line bundles $E \to X$ and $B \to X$ respectively, and a section $s\in H^0(Z, m(K_Z + E|Z)+B|Z)$ such that 
\[
\int _Z |s|^2 \omega ^{-n(m-1)} e^{-((m-1)\gamma _E + \vp _E +\vp _B)} < +\infty \quad \text{and} \quad s\wedge dT \in \si (e^{-(\vp _Z + \vp _E + \vp _B)}|Z) .
\]
Then we seek to construct a singular metric $e^{-\psi}$ such that 
\[
\int _Z |s|^2 e^{-\psi} < +\infty \quad \text{and} \quad s\wedge dT \in \si (e^{-(\vp _Z + \psi)}|Z) .
\]
If this metric satisfies condition (C), then the extension theorem will give an extension of the section $s$ to a holomorphic section $S$ of $m(K_X+Z+E)+B$ over $X$ such that $S|Z=s\wedge dT$.

Siu's idea is to use the extension theorem not only to extend $s$, but also to construct $e^{-\psi}$.  However, to get the construction off the ground, one needs to twist certain line bundles by a sufficiently positive line bundle.  Then a limiting process is used to eliminate this positive line bundle, through the method of taking powers and roots.

Our construction of the metric $e^{-\psi}$ is substantially shortened by the use of a new method of Paun introduced in \cite{p}.  Paun has eliminated the need for the use of an effective global generation of multiplier ideal sheaves, which was a long and difficult part of Siu's approach.

Finally we should mention the recent preprint \cite{c} of Claudon, who handles the case where the normal bundle of $Z$ is trivial. (Claudon has assumed that the family lies over a disk.)  Claudon's result is a direct adaptation of Siu's (or moreso Paun's) methods.  Unfortunately, such an extension theorem is of limited applicabilty when one wants to handle plruicanonical bundles by induction on dimension.  

\begin{ack}
I am grateful to Chris Hacon, who encouraged me to write this paper when he visited Stony Brook, to Mihai Paun, who graciously sent me the preprint of his simplified and generalized proof of Siu's Theorem on the deformation invariance of plurigenera, and to Mark de Cataldo, Rob Lazarsfeld, Mircea Musta\c{t}\v{a} and Mihnea Popa for helpful comments.
\end{ack}

\tableofcontents

\section{Algebro-geometric Corollaries of Theorem \ref{main}} \label{tak-section}

In this section we derive corollaries of Theorem \ref{main} that can be phrased in terms of more algebro-geometric properties of $E$ and $Z$.

\subsection{Takayama's Theorem}

As mentioned in the introduction, Theorem \ref{main} was motivated by the desire to generalize the following theorem of Takayama 

\begin{d-thm}\cite[Theorem 4.1]{tak} \label{tak-thm}
Let $X$ be a complex projective manifold, $Z \subset X$ a complex submanifold of codimension 1, and $E$ an integral divisor on $X$.  Assume that $E \sim _{\Q} A+D$ for some big and nef $Q$-divisor $A$ and some effective $\Q$-divisor $D$ such that $Z$ is in $A$-general position and $Z \not \subset {\rm Support}(D)$, and that the pair $(Z,D|Z)$ is klt.  Then the natural restriction map 
\[
H^0(X,m(K_X+Z+E)) \to H^0(Z,m(K_Z+E|Z))
\]
is surjective.
\end{d-thm}

\begin{rmk}
Recall that if $D$ is a $\Q$-divisor on $Z$, then 

\noi (i) the multiplier ideal $\si (D)$ is the multiplier ideal for the metric $e^{-\log |s_D|^2}$, where $s_D$ is the multi-section with $\Q$-divisor $D$, and 

\noi (ii) one says that the pair $(X,D)$ is {\it klt} (Kawamata Log Terminal) if $\si (D) = \co _X$.

\noi ({We say that $s$ is a multi-section of a $\Q$-divisor $D$ if there is an integer $m>0$ such that $mD$ is a $\Z$-divisor and $s^m \in H^0(X,mD)$.})
\end{rmk}
Theorem \ref{tak-thm} is a corollary of Theorem \ref{main}.  To see this, we argue as follows.  It is not hard to see (cf. Section 4 in \cite{tak}) that we may assume without loss of generality that $A$ is an ample $\Q$-divisor.  Thus we can construct a singular metric $e^{-\vp _E}$of positive curvature for $E$ as follows:  take a multiple $mA$ that is very ample, and let 
$s_D$ be the canonical multi-section of $D$ whose $\Q$-divisor is $D$.  Then we set 
\[
\vp _E = \log |s_D|^2 + \log \left ( \sum _{j=1} ^N |s_j|^{2/m} \right ),
\]
where $s_1,...,s_N$ is a basis for $H^0(X,mA)$.  Since $Z$ is not contained in the support of $D$, $e^{-\vp _E}$ restricts to $Z$ as a well-defined singular metric.  Fix any smooth metric $e^{-\vp _Z}$ for $Z$, and let $B= \co _X$ and $\vp _B \equiv 0$.  Then evidently hypotheses (P) and (T) of Theorem \ref{main} are satisfied.  Moreover, the multiplier ideal $\si (e^{-(\vp _E + \vp _Z)}|Z)$ is supported away from $Z$ because $(Z,D|Z)$ is klt.  Thus theorem \ref{tak-thm} follows.

\begin{rmk}
Under the assumption that $E$ is a big line bundle, Theorem \ref{main} can be proved by much older methods of several complex variables, going back to the work of Bombieri.  We demonstrate this in Section \ref{l2-section}.
\end{rmk}

\subsection{More general conditions on $E$}

In Theorem \ref{tak-thm}, we would like to remove the hypothesis that $E$ is big.  In some sense, this is achieved in Theorem \ref{main}.  
(Indeed, the desire to handle the case where $E$ is not necessarily big forms the initial impetus for the present article.)  However, as we mentioned, we would like to state a result that uses more intrinsic properties of the divisors, rather than a result which includes a choice of metrics.  

\begin{defn}
Let $L$ be an integral divisor on $X$.  For each integer $k > 0$, fix bases
\[
s^{(k)}_1 , ..., s^{(k)}_{N^L_ k} \in  H^0 (X, kL).
\]
Then define 
\[
\psi _L = \log \sum _{k=1} ^{\infty} \ve _k \left ( \sum _{j=1} ^{N^L_k} |s^{(k)}_j|^2 \right ) ^{2/k}.
\]
where $\ve _k > 0$ are small enough to make the sum converge.   We extend the definition to $\Q$-divisors $L$ by setting 
\[
\psi _{L} = \frac{1}{m} \psi _{mL},
\]
where $m > 0$ is the smallest integer such that $mL$ is an integral divisor.
\end{defn}

\begin{rmk}
Metrics of the form $e^{-\psi _L}$ have the following property: for every integer $m > 0$ such that $mL$ is integral, the natural inclusion 
\[
H^0(X,mL \tensor \si (e^{-m\psi _L}) ) \to H^0(X,mL)
\]
is an isomorphism.  Indeed, if for all $m>0$, $H^0(X,mL)=\{0\}$, there is nothing to prove.  On the other hand, if $\sigma = \sum c^j s_j^{(m)} \in H^0(X,mL)$ then 
\[
|\sigma |^2 e^{-m\psi _L} \le \frac{|\sigma |^2}{|s_1^{(m)}|^2 + ... + |s_{N_m}^{(m)}|^2}\le \sum |c^j|^2
\]
is bounded and thus integrable.  Thus in some sense, the metrics $e^{-\psi _L}$ have minimal singularities.
\end{rmk}

Recall that the set theoretic base locus $\text{Bs}(|L|)$ of an integral divisor $L$ is the common zero locus of all holomorphic sections of the line bundle associated to $L$.  For a more thorough discussion of this and many other matters in algebraic geometry, as well as a more algebro-geometric approach to multiplier ideals, see \cite{laz}.

\begin{d-thm}\label{tak-2}
Let $X$ be a projective algebraic manifold, $Z \subset X$ a smooth divisor, and $E$ an integral divisor on $X$ . Assume that one can write $E \sim _{\Q} E_1 + E_2$ such that the following properties hold. 
\begin{enumerate} 
\item[$(P_a)$] For some $\mu \in  \N$ such that $\mu E_1$ is integral, 
\[
{\rm Bs}(|\mu E_1|) \cup {\rm Bs}(|\mu E_1-Z |)  = \emptyset.
\] 

\item[$(T_a)$] The singular metric $e^{-\psi _{E_2}}$ restricts to a singular metric on $Z$, and 
\[
\si (e^{-\psi _{E_2}}|Z) = \co _Z.
\]
\end{enumerate} 
Then the restriction map 
\[
H^0(X, m(K_X+Z+E)) \to H^0(Z, m(K_Z + E|Z))
\]
is surjective.
\end{d-thm}

\begin{proof}[Proof of Theorem \ref{tak-2} from Theorem \ref{main}]
Let 
\[
\zeta =\psi _{\mu E_1} \quad \text{and} \quad \eta = \psi _{\mu E_1 -Z}.
\]
By hypothesis $(P_a)$, the curvature currents of the singular metrics $e^{-\zeta}$ and $e^{-\eta}$ for $\mu E_1$ and $\mu E_1-Z$ respectively are non-negative and smooth.  Let 
\[
\vp _{Z}:= \zeta - \eta, \quad \vp _{E_1} = \frac{1}{\mu} \zeta  \quad \text{and} \quad \vp _E := \vp _{E_1} + \psi _{E_2}.
\]
Note the following.
\begin{enumerate}
\item[(i)]  The metrics 
\[
\vp _{E_1}, \quad \vp _E \quad \text{and} \quad \mu \vp _{E} - \vp _Z = \mu \psi _{E_2} + \eta 
\]
have non-negative curvature currents.

\item[(ii)] $\vp _Z$ and $\psi _{E_1}$ are smooth, and thus 
\[
\si (e^{-(\vp _Z + \vp _E)}|Z) = \si (e^{-\psi _{E_2}}|Z) = \co _Z
\]
by hypothesis $(T_a)$.
\end{enumerate}
Now take $B = \co _X$, $\vp _B \equiv 0$.   The metrics $e^{-\vp _Z}$, $e^{-\vp _E}$ and $e^{-\vp _B}$ satisfy the hypotheses of Theorem \ref{main}, and moreover 
\[
\si (e^{-(\vp _Z+\vp_E+\vp _B)}|Z) = \co _Z.
\]
We thus obtain Theorem \ref{tak-2}.
\end{proof} 

\begin{rmk}
If one can write $E \sim _{\Q} E_1+E_2$ where $E_1$ is an ample $\Q$-divisor and $(Z,E_2|Z)$ is klt, then certainly properties $(P_a)$ and $(T_a)$ hold.  We thus recover Theorem \ref{tak-thm}.  
\end{rmk}

\begin{rmk}
A natural way in which the hypotheses of Theorem \ref{tak-2} might arise is the following.  Suppose $X$ and $Y$ are projective manifolds, $\text{dim} _{\C} Y > \text{dim}_{\C} X$, and $\pi : Y \to X$ is a holomorphic map whose fibers have constant dimension.  Let $Z$ and $E$ be divisors on $X$ such that $(X,Z,E)$ satisfy the hypotheses of Takayama's Theorem \ref{tak-thm}.  Then the hypotheses of Theorem \ref{tak-2} hold for $\pi ^*E$ and $\pi ^*Z$ on $Y$.
\end{rmk}

\section{$L^2$ extension}\label{l2-section}

In this section we discuss the extension of twisted canonical sections from a codimension 1 submanifold.    We begin with the local extension, which is the adjunction formula, and then pass to $L^2$ extension in the presence of certain minimal positivity hypotheses.  If the positivity hypotheses are made a little more strict, we obtain a proof of $L^2$ extension by older methods.

\subsection{Restriction of the canonical bundle}

Let $Z$ be a smooth complex hypersurface in a K\"ahler manifold $Y$, and $T \in H^0(Y,Z)$ a section whose zero divisor is $Z$.  Since the canonical bundle is the determinant of the cotangent bundle, one has 
\[
K_Y|Z = K_Z + N^*_Z,
\]
where $N^*_Z$ denotes the conormal bundle of $Z$ in $Y$, i.e., the annihilator of $T_Z$.   The adjunction formula, which amounts to saying that $dT$ is a nowhere zero section of $N^*_Z + Z$ over $Z$,  tells us that the line bundle associated to the divisor $Z$ restricts, along $Z$, to the normal bundle $N_Z$ of $Z$ in $Y$.  Equivalently,  
\[
(K_Y +Z)|Z= K_Z.
\]
More locally, if $s$ is a local section of $K_Z$, then $s\wedge dT$ is a local section of $(K_Y+Z) |Z$.

\subsection{$L^2$-extension theorem of Ohsawa-Takegoshi type}

\subsubsection{Statement of the theorem}
The setting of the theorem is the following.  Let $Y$ be a K\"ahler manifold of complex dimension $n$.  Assume there exists an analytic hypersurface $V \subset Y$ such that $Y-V$ is Stein.  Thus there are relatively compact subsets $\Omega _j \relcomp Y-V$ such that
\[
\Omega _j \relcomp \Omega _{j+1} \quad \text{and} \quad \bigcup _j \Omega _j = Y-V.
\]
Examples of such manifolds are Stein manifolds (where $V$ is empty) and projective algebraic manifolds (where one can take $V$ to be the intersection of $Y$ with a projective hyperplane in some projective space in which $Y$ is embedded).

Fix a smooth hypersurface $Z \subset Y$ such that $Z \not \subset V$.  

\begin{thm}\label{ot-thm}
Suppose given a holomorphic line bundle $H \to Y$ with a singular Hermitian metric $e^{-\kappa}$, and a singular Hermitian metric $e^{-\vp _Z}$ for the line bundle associated to the divisor $Z$, such that the following properties hold. 
\begin{enumerate}
\item[(i)] The restriction $e^{-(\kappa +\vp _Z)}|Z$ is a singular metric for $Z + H$. 
\item[(ii)] There is a global holomorphic section $T \in H^0(Y,Z)$ such that 
\[
Z= \{ T= 0\} \quad \text{and} \quad \sup_Y |T|^2e^{-\vp _Z} =1.
\]
\item[(iii)] There is an integer $\mu > 0$ such that $\mu \ii \di \dbar \kappa \ge \max (\ii \di \dbar \vp _Z, 0)$.
\end{enumerate}
Then for every $s \in H^0 (Z, K_Z + H)$ such that 
\[
\int _Z |s|^2 e^{-\kappa} < +\infty \quad \text{and} \quad s\wedge dT \in \si (e^{-(\vp _Z +\kappa)}|Z),
\]
there exists a section $S \in H^0 (Y,K_Y +Z + H)$ such that 
\[
S|Z = s \wedge dT \quad \text{and} \quad \int _Y |S|^2 e^{-(\vp _Z + \kappa)} \le 40 \pi  \mu \int _Z |s|^2 e^{-\kappa}.
\]
\end{thm}

\subsubsection{The twisted basic estimate}

Let $\Omega$ be a smoothly bounded pseudoconvex domain in our K\" ahler manifold $Y$.   Let $\tau$ and $A$ be positive smooth function on $\Omega$, and let $e^{-\psi}$ be a singular metric for $E+Z$ over  $\Omega$.  The following lemma is well-known (see, e.g., \cite{mv}).

\begin{lem}\label{mcneal-lem}
For any $(n,1)$-form $u$ in the domain of the adjoint $\dbar_\psi ^*$ of $\dbar$, the following inequality holds.
\begin{eqnarray}\label{mcneal-id}
&&\int _{\Omega} (\tau + A) \left | \dbar ^* _{\psi} u\right |^2 e^{-\psi}  + \int _{\Omega} \tau \left | \dbar u \right | ^2
e^{-\psi} \\
\nonumber && \qquad \ge \int _{\Omega} \left ( \tau \ii \di \dbar \psi -
\ii \di \dbar \tau  - \frac{1}{A} \ii \di \tau \wedge \dbar \tau \right ) (u,u)  e^{-\psi}.
\end{eqnarray}
\end{lem}

\subsubsection{Choices for $\tau$, $A$ and $\psi$, and an a priori estimate}  We follow the approach of \cite{mv}.

First, we set 
\[
\tau = a + h(a) \quad \text{and} \quad A = \frac{(1+h'(a))^2}{- h''(a)},
\]
where 
\[
h(x) = 2- x + \log (2e^{x-1} -1)
\]
and $a : \Omega \to [1,\infty)$ is a function to be chosen shortly.  Observe that for $x \ge 1$, 
\[
h'(x) = \frac{1}{2e^{x-1} -1 } \in (0,1) \quad \text{and} \quad h''(x) = \frac{- 2e^{x-1}}{(2e^{x-1}-1)^2} < 0,
\]
and thus since $1+ \log r \ge \tfrac{1}{r}$ when $r \ge 1$, 
\[
\tau \ge 1+h'(a).
\]
Moreover, $A>0$, which in necessary in our choice of $A$.  We also take this opportunity to note that 
\[
A = 2e^{a-1}.
\]
With these choices of $\tau$ and $A$, we have
\begin{eqnarray}\label{twist-curv}
- \di _{\alpha} \di _{\bar \beta} \tau - \frac{\di _{\alpha } \tau \overline{\di _{\beta} \tau}}{A} &=& - \di _{\alpha} \left ((1+h' (a)) \di _{\bar \beta} a\right ) - \frac{(1+h'(a))^2  \di _{\alpha} a \overline{\di _{\beta} a}}{A}\qquad \\
&=& (1+h'(a)) \left ( - \di _{\alpha} \di _{\bar \beta} a\right ) \nonumber 
\end{eqnarray}

Our next task is to construct the function $a$.  To this end, define
\[
v=  \log |T|^2 - \vp _Z .
\]
We note that $v \le 0$.  Fix a constant $\gamma > 1$.  We define the function $a$ to be
\[
a=a_{\ve} := \gamma -  \frac{1}{\mu}\log \left ( e^{v} + \ve ^2 \right ),
\]
where $\ve > 0$ is chosen so small that $a \ge 1$.  Later we will let $\ve$ go to $0$ and $\gamma \to 1$.

We calculate that 
\begin{eqnarray*}
-\ii \di \dbar a &=& \frac{\ii}{\mu} \di \dbar \log (e^v + \ve ^2) \\
&=& \frac{e^v}{\mu (e^v + \ve ^2)} \ii \di \dbar v + \frac{4 \ve ^2 \ii \di ( e^{v/2}) \wedge \dbar (e^{v/2})}{\mu ((e^{v/2})^2 + \ve ^2)^2}  \\
&=& -\frac{1}{\mu} \frac{e^v}{(e^v + \ve ^2)} \ii \di \dbar \vp _Z +  \frac{4\ve ^2\ii \di ( e^{v/2}) \wedge \dbar (e^{v/2})}{\mu ((e^{v/2})^2 + \ve ^2)^2}.
\end{eqnarray*}
In the last equality we have used the fact that 
\[
\ii \di \dbar v = \pi [Z]- \ii \di \dbar \vp _Z, 
\]
where $[Z]$ is the current of integration over $Z$.  The term involving the current of integration vanishes because $e^v |Z \equiv 0$.

It remains to choose the metric $e^{-\psi}$.  We take 
\[
\psi = \kappa  + \log |T|^2.
\]
Then 
\begin{eqnarray*}
&& \tau \ii \di \dbar \psi - \ii \di \dbar \tau - \frac{\ii \di \tau \wedge \dbar \tau }{A} \\
&=& \tau \ii \di \dbar \kappa + \pi[Z]\\
&& \qquad + (1+h'(a))\left (  - \frac{e^v}{\mu (e^v + \ve ^2)} \ii \di \dbar \vp _Z + \frac{4\ve ^2\ii \di ( e^{v/2}) \wedge \dbar (e^{v/2})}{\mu ((e^{v/2})^2 + \ve ^2)^2}\right )\\
&\ge & \frac{4\ve ^2\ii \di ( e^{v/2}) \wedge \dbar (e^{v/2})}{\mu ((e^{v/2})^2 + \ve ^2)^2} .
\end{eqnarray*}
The inequality follows from assumption (i) and the fact that $\tau \ge 1 +h'(a) \ge 1$.  Combining with \eqref{mcneal-id}, we obtain the following lemma.

\begin{lem}\label{basic-est-2}
Let $T = \dbar \circ \sqrt{\tau +A}$ and $S  = \sqrt{\tau} \dbar$.  Then 
for any $(n,1)$-form $u$ in the domain of the adjoint $T^*$, the following inequality holds.
\[
\int _{\Omega} \left |\left < u ,\dbar (e^{v/2})\right > \right |^2 \frac{4\ve^{2}}{\mu (e^v+\ve ^2)^{2}} e^{-\psi} \le \left ( ||T^*u||_\psi^2 + ||Su||_\psi^2 \right )
\]
\end{lem}

\subsubsection{A smooth extension and its holomorphic correction}

Since $\Omega$ is Stein, we can extend $s\wedge dT$ to a $Z+ H$-valued holomorphic $n$-form $\tilde s$ on $\Omega$.  By extending to a Stein neighborhood of $\Omega$ (which exists by hypothesis) we may also assume that 
\[
\int _{\Omega}\left | \tilde s \right |^2 e^{-(\vp _Z+\kappa)}< +\infty.
\]
(Here we have used the local integrability of $|s\wedge dT|^2 e^{-(\vp _Z + \kappa)}$ on $Z$.)  Of course, we have no better estimate on this $\tilde s$.  In particular, the estimate could degenerate as $\Omega$ grows.

In order to tame the growth of this extension $\tilde s$, we first modify it to a smooth extension.  To this end, let $\delta >0$ and let $\chi \in \cc ^{\infty} _0 ([0,1))$ be a
cutoff function with values in $[0,1]$ such that
\[
\chi \equiv 1 \ {\rm on} \ [0,\delta] \quad {\rm and} \quad |\chi
'|\le 1+\delta.
\]
We write
\[
\chi _{\ve} := \chi \left ( \frac{e^v}{\ve^2} \right ).
\]
We distinguish the smooth $(n,1)$-form with values in $Z + H$
\[
\alpha _{\ve} := \dbar \chi _{\ve} \tilde s.
\]
Then one has the estimate
\begin{eqnarray*}
\left|(u,\alpha _{\ve})_\psi\right|^2 &\le & \left ( \int _{\Omega}|\left <
u,\alpha _{\ve}\right >| e^{-\psi} \right )^2 \\
&=& \left ( \int _{\Omega} \left | \left < u, \chi ' \left ( \frac{e^v}{\ve ^2}\right ) \tilde s \wedge  \frac{2e^{v/2}\dbar (e^{v/2})}{\ve^2} \right >\right | e^{-\psi} \right )^2 \\
&\le & \mu \int _{\Omega} \left | \frac{\tilde s}{\ve ^2}
\chi ' \left ( \frac{e^v}{\ve ^2}\right )\right |^2
\frac{(e^v+\ve ^2)^2}{\ve^{2}} e^{-(\psi -v) } \\
&& \times \int _{\Omega} \left |\left < u ,\dbar (e^{v/2})\right > \right |^2 \frac{4 \ve^{2}}{\mu (e^v+\ve ^2)^{2}} e^{-\psi} \\
&\le & C_{\ve} \left ( ||T^*u||_\psi^2 + ||Su||_\psi^2 \right )
\end{eqnarray*}
where
\[
C_{\ve} := \frac{4\mu (1+\delta)^2}{\ve ^2} \int _{e^v \le \ve^2} \left | \tilde s \right |^2 e^{-(\vp _Z+\kappa)}\  {\buildrel \ve \to 0 \over \longrightarrow} \  8\pi \mu (1+\delta)^2 \int _Z |s|^2 e^{-\kappa},
\]
and the last inequality follows from Lemma \ref{basic-est-2}.

As a result of this estimate, we obtain the following theorem.

\begin{d-thm}\label{t-est-thm}
There exists a smooth $n$-form $\beta _{\ve}$ such that
\[
T \beta _{\ve} = \alpha _{\ve} \quad {\rm and} \quad \int _{\Omega}
|\beta _{\ve}|^2 e^{-\psi} \le C_{\ve}.
\]
In particular,
\[
\beta _{\ve} |Z \equiv 0.
\]
\end{d-thm}
\begin{proof} 
The inequality
\[
|(u,\alpha _{\ve})| \le C_{\ve} \left ( ||T^*u||_\psi^2 + ||Su||_\psi^2 \right )
\]
implies the continuity of the linear functional
\[
\ell : T^*u\mapsto (u,\alpha_\ve)_\psi
\] 
on $\text{Image}(T^*|\text{Kernel}(S))$.  By letting $\ell \equiv 0$ on $\text{Image}(T^*|\text{Kernel}(S))^{\perp}$, we can assume, without increasing the norm, that $\ell$ is defined on the whole Hilbert space.  The Riesz Representation Theorem then gives a solution to $T \beta _\ve = \alpha_\ve$ with the stated norm inequality on $\beta_\ve$.  The smoothness of $\beta _{\ve}$ follows from elliptic regularity.

It remains only to show that $\beta _{\ve} |Z \equiv 0$.  But one notices that $\psi$ is at least as singular as $\log |T|^2$, and thus $e^{-\psi}$ is not locally integrable at any point of $Z$.  The desired vanishing of $\beta _{\ve}$ follows.
\end{proof}

\begin{proof}[Conclusion of the proof of Theorem \ref{ot-thm}]
We note first that 
Now, 
\begin{eqnarray*}
\lim _{\ve \to 0}  \sup _{\Omega} (e^{v} (\tau +A)) &=& \lim _{\ve \to 0} \sup _{\Omega} e^{\gamma -a} ( 2 + \log (2e^{a-1} -1) + 2e^{a-1})\\
&\le & \sup _{x \ge 1} e^{\gamma -x} ( 3+2 e^{x-1}) \\
&=& 5e^{\gamma -1}. 
\end{eqnarray*}
Next let 
\[
S_{\ve} := \chi _{\ve} \tilde s - \sqrt{\tau + A} \beta _{\ve}.
\]
Then $S_{\ve}$ is a holomorphic section, $S_{\ve} |Z = \tilde s |Z = s\wedge dT$, and we have the estimate
\[
\int _{\Omega} |S_{\ve}|^2 e^{-(\vp _Z + \kappa)} \le 5e^{\gamma -1} (1+ o(1))8\pi \mu (1+\delta )^2 \int _Z |s|^2 e^{-\kappa} \quad \text{as } \ve \to 0.
\]
Here we have used the estimate obtained in Theorem \ref{t-est-thm}.  
The term $o(1)$ comes in because $\chi _{\ve} \tilde s$ is a smooth, uniformly bounded section whose support approaches a set of measure zero, and thus the integral of the square norm of $\chi _{\ve} \tilde s$ converges to zero.

Now, by the sub-mean value property, uniform $L^2$-esimates with plurisubharmonic (hence locally bounded above) weights implies locally uniform sup-norm estimates.  It follows from Montel's Theorem that $S_{\ve}$ converges to a holomorphic section $S$.  Evidently 
\[
S|Z\cap \Omega = s\wedge dT
\]
and
\[
\int _{\Omega} |S|^2 e^{-(\vp _Z + \kappa)} \le 5e^{\gamma -1} \times8\pi \mu (1+\delta )^2 \int _Z |s|^2 e^{-\kappa}.
\]
Finally, by the uniformity of all estimates, we may let $\delta \to 0$, $\gamma \to 1$ and then $\Omega \to Y$.  The proof is complete.
\end{proof}

\subsection{$L^2$-extension of H\" ormander-Bombieri-Skoda type}

Let $Y$ be a compact complex algebraic manifold and $Z\subset Y$ a smooth complex hypersurface. 

\begin{thm}\label{hbs-thm}
Suppose given a holomorphic line bundle $H \to Y$ with a singular Hermitian metric $e^{-\kappa}$, and a smooth Hermitian metric $e^{-\vp _Z}$ for the line bundle associated to the divisor $Z$, such that the following properties hold. 
\begin{enumerate}
\item[(i)] The restriction $e^{-(\kappa +\vp _Z)}|Z$ is a singular metric for $Z + H$. 
\item[(ii)] There is a global holomorphic section $T \in H^0(Y,Z)$ such that 
\[
\sup _Y |T|^2e^{-\vp _Z} =1.
\]
\item[(iii)] There is a constant $c > 0$ such that $\ii \di \dbar \kappa \ge c$.
\end{enumerate}
Then there exists a constant $C$ depending $Z$ and $Y$ but not on $H$, with the following property.  For every $s \in H^0 (Z, K_Z + H)$ such that 
\[
\int _Z |s|^2 e^{-\kappa} < +\infty \quad \text{and} \quad s\wedge dT \in \si (e^{-(\vp _Z +\kappa)})\cdot \co _Z
\]
there exists a section $S \in H^0 (Y,K_Y +Z + H)$ such that 
\[
S|Z = s \wedge dT \quad \text{and} \quad \int _Y |S|^2 e^{-(\vp _Z + \kappa)} \le \frac{C}{c} \int _Z |s|^2 e^{-\kappa}.
\]
\end{thm}

The approach we take to prove Theorem \ref{hbs-thm}  is as follows.  We first extend the section $s$ locally, with uniform estimates.  The local extensions need not agree, so we correct their discrepancy by solving a Cousin I problem with good $L^2$ bounds.

\subsubsection{Local extension}

\begin{lem}\label{loc-extension-lemma}
Let $U$ be a smoothly bounded open set in $Y$ such that 
\begin{enumerate}
\item[(a)] $U$ is biholomorphic to the unit ball, and 
\item[(b)] $Z \cap U$ is the hyperplane $z^n = 0$ in the coordinates on the ball $U$.
\end{enumerate}
Let $U' \relcomp U$ be the ball of radius $1/2$.  Then for every section $s \in H^0(U\cap Z, K_Z +H)$ such that 
\[
\int _{Z\cap U} |s|^2 e^{-\kappa} < +\infty \quad \text{and} \quad s\wedge dT \in \si (e^{-(\vp _Z +\kappa)}|{Z\cap U'}),
\]
there exists a section $\tilde s \in H^0(U' , K_X+Z+H)$ such that 
\[
\tilde s |(U'\cap Z) = s_o \wedge dT \quad \text{and} \quad \int _{U'} |\tilde s|^2 e^{-(\vp _Z + \kappa)} \le C_o \int _{Z\cap U'} |s|^2 e^{-\kappa}.
\]
The constant $C_o$ is independent of $s$.
\end{lem}

One could deduce a stronger result directly from the classical Ohsawa-Takegoshi extension theorem in \cite{ot-87}, but Lemma \ref{loc-extension-lemma} is far more elementary.  We leave the proof as an exercise to the reader.

\subsubsection{An $L^2$ Cousin I problem}

Let us cover $Y$ by a finite number of coordinate unit balls $U_j$, such that the concentric balls $V_j$ of radius $1/2$ also cover $Y$.   Suppose we have sections $\tilde s_j \in H^0(V_j , K_X +Z+H)$ such that 
\[
\tilde s_j |V_j \cap Z = s \wedge dT |V_j \cap Z.
\]
(If $V_j \cap Z = \emptyset$, then the restriction condition is vacuously satisfied.)  Consider the cocycle 
\[
G_{ij} = \tilde s_i - \tilde s_j \quad \text{on  }U_i \cap U_j.
\]
Note that 
\[
\int _{V_i \cap V_j } |G_{ij}|^2 e^{-(\vp _Z + \kappa)} \le C \int _{V_j \cup V_j\cap Z} |s|^2 e^{-\kappa} \quad \text{and} \quad  G_{ij}|U_i \cap U_{j}\cap Z = 0.
\]
Thus, since $dT |Z \neq 0$, we see that 
\[
G_{ij} / T = f_{ij} \in H^0(U_{i} \cap U_{j} , K_Y +H)
\]
and 
\[
\int _{V_i \cap V_j} |G_{ij}|^2 e^{-(\kappa+\log |T|^2)} = \int _{V_i \cap V_j}  |f_{ij}|^2 e^{-\kappa}
\]

We seek holomorphic sections $g_j \in H^0 (V_j, K_Y +Z+H)$ such that 
\[
\int _{V_j }|g_j|^2 e^{-(\vp _Z+\kappa)} \le K_1 \int _{Z}|s|^2 e^{-(\vp _Z + \kappa)}, 
\]
\[
g_i - g_j = G_{ij} \quad \text{and} \quad g_i |V_j \cap Z = 0.
\]
If such sections are found, then the section $s \in H^0 (Y,K_Y +Z+H)$ defined by 
\[
s = s_i -g_i \quad \text{on} \ V_i
\]
gives the proof of Theorem \ref{hbs-thm}.

\subsubsection{Solution of the Cousin I problem}  We will use the convention that the constant $C$ may change from line to line.  We leave it to the reader to check that all uses of constants in the estimates do not depend on the section $s$, but only on properties of $Y$ and $Z$.

Let $\{ \chi _j \}$ be a partition of unity subordinate to the open cover $\{V_j \}$.   Consider the smooth sections 
\[
h_j = \sum _{k} G_{jk} \chi _{k}.
\]
Then 
\[
h_i - h_j = \sum _k (G_{ik} +G_{kj}) \chi _k = \sum _k G_{ij} \chi _k = G_{ij}.
\]
It follows that the differential $(0,1)$-form $\alpha$ with values in $K_Y +Z+H$ defined by
\[ 
\alpha = \dbar h_j \quad \text{on  }U_j
\]
is well-defined.  Moreover
\[
\int _Y |\alpha |^2 e^{-(\kappa + \log |T|^2)} \le C \sum _{i,j} \int _{V_i \cap V_j} |f_{ij}|^2 e^{-\kappa}
\]
Recalling that $V_j$ is the ball of radius $1/2$, let $W_j \subset V_j$ denote the spherical shells of inner radius $1/4$ and outer radius $1/2$.  Then 
\begin{eqnarray*}
\int _{V_i \cap V_j} |f_{ij}|^2 e^{-\kappa} &\le & C \int _{W_i \cap W_j} |f_{ij}|^2 e^{-\kappa} \\
&\le & C \int _{W_j \cap W_j} |T|^2 |f_{ij}|^2 e^{-(\vp _Z+\kappa)}\\
&\le & C\int _{V_i \cap V_j} |G_{ij}|^2 e^{-(\vp _Z+\kappa)}\\
&\le & C\int _{V_i \cup V_j} |s|^2 e^{-\kappa}.
\end{eqnarray*}
(The first inequality follows from an appropriate application of the sub-mean value property.)  Thus 
\[
\int _Y |\alpha |^2 e^{-(\kappa + \log |T|^2)} \le C \int _{Z} |s|^2 e^{-\kappa}.
\]

By H\" ormander's Theorem, there is a smooth section $u$ such that 
\[
\dbar u = \alpha \quad \text{and} \quad \int _Y |u|^2 e^{-(\vp _Z+\kappa)} \le \int _Y |u|^2 e^{-(\kappa +\log |T|^2)} \le \frac{C}{c} \int _Z|s|^2 e^{-\kappa}, 
\]
where $\ii \di \dbar \kappa > c \omega$, and in the first inequality we have used $|T|^2e^{-\vp _Z} \le 1$.  Observe that the second inequality forces $u$ to vanish along $Z$.

We let $g_i = h_i - u$.  Then $g_i$ are holomorphic, vanish along $V_j \cap Z$, and satisfy the estimate
\[
\int _{V_i}|g_i|^2 e^{-(\vp _Z+\kappa)} \le \frac{C}{c} \int _Z |s|^2 e^{-(\vp _Z + \kappa)} .
\]
The proof of Theorem \ref{hbs-thm} is complete.\qed

\section{Proof of Theorem \ref{main}}\label{last-section}

For the rest of the paper, we normalize our canonical section $T$ of the line bundle associated to $Z$, so that 
\[
\sup _X |T|^2 e^{-\vp _Z} = 1.
\]

\begin{rmk}
Let $s \in H^0(Z, \co _Z(m(K_Z + E|Z)+B|Z)\tensor \si (e^{-(\vp _B+\vp _B+\vp _Z)}|Z))$ be the section to be extended.  We note that in fact, 
\[
\int _Z |s|^2 \omega ^{-(n-1)(m-1)}e^{-(m-1)\gamma _E} e^{-(\vp _E + \vp _B)} < +\infty, 
\]
since 
\[
e^{-(\vp _E+\vp _B)} \le \left ( \sup _Z e^{\vp _Z - \gamma _Z} \right ) e^{\gamma _Z} e^{-(\vp _E+\vp _B+\vp _Z)} \le C  e^{\gamma _Z} e^{-(\vp _E+\vp _B+\vp _Z)}.
\]
The last inequality follows since, by our convention, the local potentials of singular metrics are locally uniformly bounded above.
\end{rmk}

\subsection{Paun's induction}

Fix a holomorphic line bundle $A \to X$ sufficiently positive as to have the following property:
\begin{enumerate}
\item[(GG)]  For each $0 \le p \le m-1$ the global sections $H^0(X, p(K_{X} + Z + E)+A)$ generate the sheaf $\co _X(p(K_{X} + Z + E)+A)$.
\end{enumerate}

\noi Let us fix bases 
\[
\{\tilde \sigma ^{(p)}_j \ ;\ 1 \le j \le N_p \}
\]
of $H^0(X, p(K_{X}+Z+E)+ A )$.  We let $\sigma ^{(p)} _j \in H^0(Z,p(K_{Z} + E|Z)+ A|Z)$ be such that 
\[
\tilde \sigma ^{(p)} _j |Z = \sigma ^{(p)} _j \wedge (dT) ^{\tensor p}.
\]
We also fix smooth metrics  
\[
e^{-\gamma _Z},\ e^{-\gamma _E} \ \text{and } e^{-\gamma _B}\text{ for } Z \to X,\ E\to X \text{ and }B \to X
\]
respectively.

\begin{prop}\label{paun-prop}
There exist a constant $C < +\infty$ and sections 
\[
\{ \tilde \sigma ^{(km+p)}_j \in H^0(X,(km+p) (K_X+Z+E)+kB+A) \ ;\ 1 \le j \le N_p\}_{0 \le p \le m-1 , k = 0 , 1 , 2 ,...}
\]
with the following properties.  
\begin{enumerate}
\item[(a)] $\tilde \sigma _j ^{(mk+p)} | Z = s^{\tensor k} \tensor \sigma _j ^{(p)} \wedge (dT) ^{(km+p)}$

\item[(b)]  If $k \ge 1$, 
\[
\int _X \frac{\sum _{j=1} ^{N_0} |\tilde \sigma _j ^{(mk)}|^2 e^{-(\gamma _Z + \gamma _E +\gamma _B)}}{\sum _{j=1} ^{N_{m-1}} |\tilde \sigma _j ^{(mk-1)}|^2}  \le C.
\]

\item[(c)]  For $1 \le p \le m-1$,
\[
\int _X \frac{\sum _{j=1} ^{N_p} |\tilde \sigma _j ^{(mk + p)}|^2e^{-(\gamma _Z + \gamma _E)}}{\sum _{j=1} ^{N_{p-1}} |\tilde \sigma _j ^{(mk + p-1)}|^2}  \le C.
\]

\end{enumerate}
\end{prop}

\begin{proof} (Double induction on $k$ and $p$.)  Fix a constant $\widehat C$ such that the 
\[
\sup _{X} \frac{\sum _{j=1}^{N_{0}} |\tilde \sigma ^{(0)} _j |^2\omega ^{n(m-1)}e^{(m-1)(\gamma _Z+\gamma _E)}}{\sum _{j=1}^{N_{m-1}} |\tilde \sigma ^{(m-1)} _j |^2} \le \widehat C
\]
and 
\[
\sup _{Z} \frac{\sum _{j=1}^{N_{0}} |\sigma ^{(0)} _j |^2\omega ^{(n-1)(m-1)}e^{(m-1)\gamma _E}}{\sum _{j=1}^{N_{m-1}} |\sigma ^{(m-1)} _j |^2} \le \widehat C, 
\]
and  for all $0 \le p \le m-2$, 
\[
\sup _X  \frac{\sum _{j=1}^{N_{p+1}} |\tilde \sigma ^{(p+1)} _j |^2\omega ^{-n}e^{-(\gamma _Z+\gamma _E)}}{\sum _{j=1}^{N_p} |\tilde \sigma ^{(p)} _j |^2} \le \widehat C, 
\]
and
\[
\sup _Z  \frac{\sum _{j=1}^{N_{p+1}} |\sigma ^{(p+1)} _j |^2\omega ^{-(n-1)}e^{-\gamma _E}}{\sum _{j=1}^{N_p} |\sigma ^{(p)} _j |^2} \le \widehat C.
\]

\noi ($k=0$)
As far as extension there is nothing to prove.  Note that
\[
\int _X \frac{\sum _{j=1} ^{N_p} |\tilde \sigma _j ^{(p)}|^2e^{-(\gamma _Z + \gamma _E)}}{\sum _{j=1} ^{N_{p-1}} |\tilde \sigma _j ^{(p-1)}|^2}  \le \widehat C \int _X \omega ^n.
\]

\noi ($k \ge 1$)  Assume the result has been proved for $k-1$.

\noi \underline{(($p=0$))}:   Consider the sections $s^{\tensor k} \tensor \sigma _j ^{(0)}$, and define the semi-positively curved metric 
\[
\psi _{k,0} := \log \sum _{j=1} ^{N_{m-1}} |\tilde \sigma ^{(km-1)}_j|^2 
\]
for the line bundle $(mk-1)(K_X+Z+E)+(k-1)B+A$.  Observe that locally, 
\begin{eqnarray*}
|(s\wedge dT ^m)^k\tensor \sigma _j ^{(0)}|^2 e^{-(\vp _Z+\psi_{k,0}+\vp _E + \vp _B)}&=& |s \wedge dT ^m|^2 \frac{|\sigma ^{(0)} _j|^2 e^{-(\vp_Z + \vp _E + \vp _B)}}{\sum _{j=1} ^{N_{m-1}} |\sigma ^{(m-1)}_j|^2 } \\
&\lesssim& |s|^2e^{-(\vp_Z + \vp _E + \vp _B)}.
\end{eqnarray*}
Moreover, we have 
\[
\mu \ii \di \dbar (\psi _{k,0} + \vp _E + \vp _B) \ge \max (\ii \di \dbar \vp _Z , 0).
\]
Finally,
\begin{eqnarray*}
&& \int _{Z} |s^k\tensor \sigma _j ^{(0)}|^2 e^{-(\psi_{k,0}+\vp _E + \vp _B)} \\
&=& \int _{Z} |s|^2 \frac{|\sigma ^{(0)} _j|^2 e^{(m-1) \gamma _E} e^{-((m-1) \gamma _E + \vp _E + \vp _B)}}{\sum _{j=1} ^{N_{m-1}} |\sigma ^{(m-1)}_j|^2 }< +\infty.
\end{eqnarray*}
In view of the remark at the beginning of Section \ref{last-section}, we may thus apply Theorem \ref{ot-thm} to obtain sections 
\[
\tilde \sigma ^{(km)}_j \in H^0(X, mk (K_X +Z+E)+kB+A), \quad 1 \le j \le N_0
\]
such that 
\[
\tilde \sigma ^{(km)}_j  |Z = s^{\tensor k} \tensor \sigma ^{(0)}_j \wedge (dT)^{\tensor km} , \quad 1 \le j \le N_0,
\]
and 
\[
\int _{X} |\tilde \sigma ^{(km)}_j|^2 e^{-(\psi_{k,0}+\vp_Z + \vp_E + \vp_B)} \le 40\pi \mu \int _{Z} |s|^2 \frac{|\sigma ^{(0)} _j|^2 e^{-(\vp _E + \vp _B)}}{\sum _{j=1} ^{N_{m-1}} |\sigma ^{(m-1)}_j|^2 }.
\]
Summing over $j$, we obtain
\begin{eqnarray*}
&& \int _{X} \frac{\sum _{j=1} ^{N_o} |\tilde \sigma ^{(km)}_j|^2 e^{-(\gamma_Z + \gamma _E + \gamma _B)}} {\sum _{j=1} ^{N_{m-1}} |\tilde \sigma ^{(km-1)}_j|^2} \\
&\le& \sup _X e^{\vp _Z + \vp _E + \vp _B - \gamma_Z - \gamma _E -\gamma _B} \int _{X} \frac{\sum _{j=1} ^{N_o} |\tilde \sigma ^{(km)}_j|^2 e^{-(\vp _Z + \vp _E + \vp _B)}} {\sum _{j=1} ^{N_{m-1}} |\tilde \sigma ^{(km-1)}_j|^2}\\
&\le&40 \pi \sup _X e^{\vp _Z + \vp _E + \vp _B - \gamma_Z - \gamma _E -\gamma _B} \int _{Z} |s|^2 \frac{\sum _{j=1} ^{N_0} |\sigma ^{(0)} _j|^2e^{-(\vp _E + \vp _B)}}{\sum _{j=1} ^{N_{m-1}} |\sigma ^{(m-1)}_j|^2 } e^{-\kappa}\\
&\le&40\pi \widehat C \sup _X e^{\vp _Z + \vp _E + \vp _B - \gamma_Z - \gamma _E -\gamma _B} \int _{Z} |s|^2 \omega ^{-(n-1)(m-1)}e^{-((m-1) \gamma _E + \vp _E + \vp _B)}
\end{eqnarray*} 

\noi \underline{(($1 \le p\le m-1$))}:  Assume that we have obtained the sections $\tilde \sigma ^{(km+p-1)}_j, \ 1\le j \le N_{p-1}$.  Consider the non-negatively curved singular metric 
\[
\psi _{k,p} := \log \sum _{j=1} ^{N_{p-1}} |\tilde \sigma ^{(mk+p-1)}_j |^2
\]
for $(km +p-1) (K_X +Z+E)+kB +A$.  We have 
\[
|s^k\tensor \sigma _j ^{(p)}|^2 e^{-(\vp _Z +\psi_{k,p} + \vp _E)} = \frac{|\sigma ^{(p)} _j|^2e^{-(\vp _Z+\vp _E)}}{\sum _{j=1} ^{N_{p-1}} |\sigma ^{(p-1)}_j|^2 } \lesssim e^{-(\vp _Z+\vp _E)},
\]
which is locally integrable by the hypothesis (T).  Next, 
\begin{eqnarray*}
\int _{Z} |s^k\tensor \sigma _j ^{(p)}|^2 e^{-(\psi_{k,p} + \vp _E)} &=&  \int _{Z} \frac{|\sigma ^{(p)} _j|^2e^{-\vp _E}}{\sum _{j=1} ^{N_{p-1}} |\sigma ^{(p-1)}_j|^2 }\\
&\le & C^{\star} \int _Z  e^{\gamma _Z} \frac{|\sigma ^{(p)} _j|^2e^{-(\vp _Z+\vp _E)}}{\sum _{j=1} ^{N_{p-1}} |\sigma ^{(p-1)}_j|^2}< +\infty,
\end{eqnarray*}
where 
\[
C^{\star} := \sup _Z e^{\vp _Z - \gamma _Z}.
\]
By Theorem \ref{ot-thm} there exist sections 
\[
\tilde \sigma ^{(km+p)}_j \in H^0(X, (mk+p)(K_X +Z+E)+kB +A), \quad 1 \le j \le N_0
\]
such that 
\[
\tilde \sigma ^{(km+p)}_j  |Z = s^{\tensor k} \tensor \sigma ^{(p)}_j \wedge (dT)^{\tensor km +p} , \quad 1 \le j \le N_p,
\]
and 
\[
\int _{X} |\tilde \sigma ^{(km+p)}_j|^2 e^{-(\psi_{k,p}+\vp _Z +\vp _E)} \le 40 \pi \mu \int _{Z}\frac{|\sigma ^{(p)} _j|^2e^{-\vp _E}}{\sum _{j=1} ^{N_{p-1}} |\sigma ^{(p-1)}_j|^2 }.
\]
Summing over $j$, we obtain
\[
\int _{X} \frac{\sum _{j=1} ^{N_p} |\tilde \sigma ^{(km+p)}_j|^2e^{-(\gamma _Z +\gamma _E)}} {\sum _{j=1} ^{N_{p-1}} |\tilde \sigma ^{(km+p-1)}_j|^2} \le 40\pi \mu \sup _X e^{\vp _Z + \vp _E - \gamma _Z - \gamma _E} \widehat C \int _{Z}e^{-\vp _E} \omega ^{n-1}.
\]
Letting $C$ be the maximum of the numbers 
\begin{eqnarray*}
&& \widehat C \int _X \omega ^n,\\
&& 40\pi \widehat C \sup _X e^{\vp _Z + \vp _E + \vp _B - \gamma_Z - \gamma _E -\gamma _B} \int _{Z} |s|^2 \omega ^{-{n-1}(m-1)}e^{-((m-1) \gamma _E + \vp _E + \vp _B)}\\
\text{and}\quad &&  40\pi \mu \sup _X e^{\vp _Z + \vp _E - \gamma _Z - \gamma _E} \widehat C \int _{Z}e^{-\vp _E} \omega ^{n-1}
\end{eqnarray*}
completes the proof.
\end{proof}

\subsection{Siu's Construction of the metric}

Fix a smooth metric $e^{-\psi}$ for $A \to X$.  Consider the functions 
\[
\lambda _{N} := \log \sum _{j=1} ^{N_p} |\tilde \sigma ^{(km+p)} _j |^2\omega ^{-n(mk+p)} e^{-(km (\gamma _Z + \gamma_E) +k \gamma _B + \psi)},
\]
where $N = km + p$.  We then have the following lemma.

\begin{lem}\label{jensen}
For any non-empty open subset $V \subset X$ and any smooth function $f: \overline{V} \to \R _+$,
\begin{eqnarray*}
\frac{1}{\int _V f \omega ^{n}} \int _V (\lambda _N - \lambda _{N-1}) f \omega ^{n} \le \log \left ( \frac{C\sup _V f}{\int _V f \omega ^{n}} \right ).
\end{eqnarray*}
\end{lem}

\begin{proof}
Observe that by Proposition \ref{paun-prop}, there exists a constant $C$ such that for any open subset $V \subset X$,
\[
\int _V (e^{\lambda _N - \lambda _{N-1}}) f \omega ^{n} \le C \sup _V f.
\]
The lemma follows from an application of (the concave version of) Jensen's inequality to the concave function $\log$.
\end{proof}

Consider the function 
\[
\Lambda _k = \frac{1}{k} \lambda _{mk}.
\]
Note that $\Lambda _k$ is locally the sum of a plurisubharmonic function and a smooth function.  By applying Lemma \ref{jensen} and using the telescoping property, we see that for any open set $V \subset X$ and any smooth function $f : \overline{V} \to \R _+$, 
\begin{equation}\label{L-bound}
\frac{1}{\int _V f \omega ^n } \int _V \Lambda _k f \omega ^n\le m \log \left ( \frac{C\sup _V f}{\int _V f\omega ^n } \right ) .
\end{equation}

\begin{prop}
There exists a constant $C_o$ such that 
\[
\Lambda _k (x) \le C_o, \quad x \in X.
\]
\end{prop}

\begin{proof}
Let us cover $X$ by coordinate charts $V_1,...,V_N$ such that for each $j$ there is a biholomorphic map $F_j$ from $V_j$ to the ball $B(0,2)$ of radius $2$ centered at the origin in $\C ^n$, and such that if $U_j = F^{-1}_j(B(0,1))$, then $U_1,...,U_N$ is also an open cover.  Let $W_j = V_j \setminus F^{-1}_j  (B(0,3/2))$.

Now, on each $V_j$, $\Lambda _k$ is the sum of a plurisubharmonic function and a smooth function.  Say $\Lambda _k = h + g$ on $V_j$, where $h$ is plurisubharmonic and $g$ is smooth.  Then for constant $A_j$ we have 
\begin{eqnarray*}
\sup _{U_j} \Lambda _k &\le& \sup _{U_j} g  + \sup _{U_j} h \\
&\le& \sup _{U_j} g  + A_j  \int_{W_j}h \cdot  F_{j*} dV\\
&\le& \sup _{U_j} g - A_j \int_{W_j}g \cdot  F_{j*} dV + A_j \int_{W_j}\Lambda _k \cdot  F_{j*} dV\\
\end{eqnarray*}
Let 
\[
C_j := \sup _{U_j} g - A_j  \int_{W_j}g \cdot  F_{j*} dV
\]
and define the smooth function $f_j$ by 
\[
f_j \omega ^n = F_{j*} dV.
\]
Then by \eqref{L-bound} applied with $V= W_j$ and $f=f_j$, we have
\[
\sup _{U_j} \Lambda _k \le C_j + m A_j \log \left ( \frac{C\sup _{W_j} f_j}{\int _{W_j} f_j \omega ^n} \right ) \int _{W_j} f_j \omega ^n.
\]
Letting 
\[
C_o := \max _{1\le j \le N} \left \{ C_j + m A_j \log \left ( \frac{C\sup _{W_j} f_j}{\int _{W_j} f_j \omega ^n} \right ) \int _{W_j} f_j \omega ^n \right \}
\]
completes the proof.
\end{proof}

Since the upper regularization of the lim sup of a uniformly bounded sequence of plurisubharmonic functions is plurisubharmonic (see, e.g., \cite[Theorem 1.6.2]{h}), we essentially have the following corollary.

\begin{cor}
The function 
\[
\Lambda (x) := \limsup _{y \to x} \limsup _{k \to \infty} \Lambda _k(y)
\]
is locally the sum of a plurisubharmonic function and a smooth function.  
\end{cor}

\begin{proof}
One need only observe that the function $\Lambda _k$ is obtained from a singular metric on the line bundle $m(K_X+Z+E)+B$  (this singular metric $e^{-\kappa _k}$ will be described shortly)  by multiplying by a fixed smooth metric of the dual line bundle. 
\end{proof}

Consider the singular Hermitian metric $e^{-\kappa}$ for $m(K_X+Z+E)+B$ defined by
\[
e^{-\kappa} = e^{-\Lambda} \omega ^{-nm}e^{-(m(\gamma _Z+\gamma _E) + \gamma _B)}.
\]
This singular metric is given by the formula
\[
e^{-\kappa (x)} = \exp \left (- \limsup _{y\to x} \limsup _{k \to \infty} \kappa _k (y) \right ),
\]
where
\[
e^{-\kappa _k} = e^{-\Lambda_k} \omega ^{-nm}e^{-(m(\gamma_Z+\gamma _E)+\gamma _B)}.
\]
The curvature of $e^{-\kappa _k}$ is thus
\begin{eqnarray*}
\ii \di \dbar \kappa _k &=& \frac{\ii}{k} \di \dbar \log \sum _{j=1} ^{N_0} |\tilde \sigma ^{(mk)}_j |^2 - \frac{1}{k} \ii \di \dbar \psi \\
&\ge &  - \frac{1}{k} \ii \di \dbar \psi 
\end{eqnarray*}

We claim next that the curvature of $e^{-\kappa}$ is non-negative.  To see this, it suffices to work locally.  Then we have that the functions 
\[
\kappa _k + \frac{1}{k} \psi
\]
are plurisubharmonic.  But
\[
 \limsup _{y\to x} \limsup _{k \to \infty}\kappa _k + \frac{1}{k} \psi =  \limsup _{y\to x} \limsup _{k \to \infty} \kappa _k = \kappa.
 \]
It follows that $\kappa$ is plurisubharmonic, as desired.

\subsection{Conclusion of the proof}
Notice that, after identifying $K_Z$ with $(K_X+Z)|Z$ by dividing by $dT$, 
\[
\kappa _k |Z = \log |s|^2 + \frac{1}{k} \log \sum _{j=1} ^{N_0} |\sigma ^{(0)} _j |^2.
\]
Thus we obtain
\[
e^{-\kappa} |Z = \frac{1}{|s|^2}.
\]
It follows that 
\begin{eqnarray*}
\int _{Z} |s|^2 e^{-\tfrac{(m-1)\kappa + \vp _E + \vp _B}{m}} &=& \int _{Z} |s|^{2/m} e^{- \frac{\vp _E+\vp _B}{m}} \\
&\le&  \left ( \int _{Z} \omega ^{n-1} \right ) ^{\tfrac{m-1}{m}} \!\!\!\!\left ( \int _{Z} |s|^2 \omega ^{-(n-1)(m-1)} e^{-(\vp _E+\vp _B)} \right ) ^{\tfrac{1}{m}} \\
&<& +\infty,
\end{eqnarray*}
where the first inequality is a consequence of H\" older's Inequality.   Next, on $Z$, 
\[
 |s \wedge dT|^2 e^{-(\vp _Z+ \tfrac{(m-1)\kappa + \vp _E + \vp _B}{m})} \sim |s|^{2/m} e^{- \frac{\vp _Z+\vp _E+\vp _B}{m}} e^{-\tfrac{m-1}{m} \vp _Z}.
\]
Now, by another application of H\"older's Inequality,  we have (locally on $Z$) that 
\[
\int |s|^{2/m} e^{- \frac{\vp _Z+\vp _E+\vp _B}{m}} e^{-\tfrac{m-1}{m} \vp _Z} \le \left ( \int |s|^2 e^{-(\vp _Z + \vp _E + \vp _B)} \right ) ^{1/m} \times \left ( \int e^{-\vp _Z}\right ) ^{(m-1)/m}.
\]
Since $e^{-\vp _Z} \le C e^{-(\vp _E+\vp _Z)}$ is locally integrable, we obtain the local integrability of 
\[
|s|^{2/m} e^{- \frac{\vp _Z+\vp _E+\vp _B}{m}} e^{-\tfrac{m-1}{m} \vp _Z}.
\]
Finally,
\[
\mu m \ii \di \dbar \left ( \tfrac{m-1}{m} \kappa + \tfrac{\vp _E+\vp _B}{m}\right ) \ge \mu \ii \di \dbar (\vp _E+\vp _B) \ge \max (\ii \di \dbar \vp _Z ,0).
\]
An application of Theorem \ref{ot-thm} completes the proof of Theorem \ref{main}.\qed

\end{document}